\documentclass[a4paper]{amsart}
\makeatother
\usepackage{amsmath}
\usepackage{amsthm}
\theoremstyle{definition}
\usepackage[english]{babel}
\usepackage{amsmath}
\usepackage{tikz-cd}
\usepackage{amsfonts}
\usepackage{amssymb}
\usepackage{algorithm2e}
\usepackage{array}
\usepackage{graphicx}
\usepackage{amsthm,amssymb}
\usepackage{epsfig}
\usepackage{float}
\usepackage{color}
\usepackage{hyperref}
\hypersetup{
    colorlinks=true,
    linkcolor=black,
    filecolor=blue,      
    urlcolor=blue,
    }
\usepackage[top=1in, bottom=1.4in, left=1.4in, right=1.4in]{geometry}
\usepackage{subfigure}
\numberwithin{figure}{section}
\newtheorem{theorem}{Theorem}[section]
\newtheorem*{theorem*}{Theorem \ref{XY}}
\newtheorem*{corollary*}{Corollary \ref{XZ}}
\newtheorem*{theorem**}{Theorem \ref{QQ}}
\newtheorem{lemma}[theorem]{Lemma}

\newtheorem{corollary}[theorem]{Corollary}
\newtheorem{proposition}[theorem]{Proposition}
\newtheorem{example}{Example}[section]
\newtheorem{claim}{Claim}
\newtheorem{remark}[theorem]{Remark}

\newtheorem{question}{Question}[section]
\newtheorem*{question***}{Question \ref{DDD}}
\date{}
\author{}
\usepackage{tikz}
\usetikzlibrary{matrix,arrows,shapes.geometric}
\usepackage{longtable}
\usepackage{soul}
\usepackage{enumerate}
\newpage
\title{Buchsbaumness, Macaulayfication and Castelnuovo-Mumford regularity of  monomial curves}
 \author{Biplab Dawn and Kumari Saloni }
\date{}

\begin{document}

\address{Indian Institute of Technology Patna, Bihar, India 801106}
\email{biplab\_2321ma10@iitp.ac.in}
\address{Indian Institute of Technology Patna, Bihar, India 801106}
\email{ksaloni@iitp.ac.in}
\date{}
\thanks{ The second author is supported by ANRF Power grant No.: $\mbox{SPG}/2022/002099$, Govt. of India. }
\thanks{Both the authors declare equal contribution in the paper.}
\dedicatory{}
\begin{abstract}
Projective monomial curves are associated with rings generated by monomials of equal degree in two variables. In this paper, we give a  class of non-smooth, non Cohen-Macaulay $k$-Buchsbaum projective monomial curves for any $k\geq 1$ and find the monomial generators for the respective  Macaulayfication. Moreover, we demonstrate a general method to find the Macaulayfication of a $k$-Buchsbaum monomial curve for any $k\geq 1$. We also discuss Castelnuovo-Mumford regularity of certain curves in terms of $k$-Buchsbaumness and reduction number.    
\end{abstract}
\maketitle
\section{Introduction}
Let $\mathcal{K}$ be a field. A curve $\mathcal{C}$ is a one-dimensional closed subscheme of $\mathbb{P}^{p+1}$ that is locally Cohen-Macaulay and equidimensional. Let  
\[
M(\mathcal{C}) = \bigoplus_{m \in \mathbb{Z}} H^1 (\mathbb{P}^{p+1}_\mathcal{K}, I_{\mathcal{C}}(m))
\]  
be the Hartshorne-Rao module of $\mathcal{C}$, where $I_{\mathcal{C}}$ denotes the ideal sheaf of $\mathcal{C}$. It is known that $M(\mathcal{C}) \cong H^1_\mathfrak{M}(R)$ (see \cite{8} and \cite[Corollary 4.7]{3}), where $R = \mathcal{K}[x_0, \dots, x_{p+1}]/I(\mathcal{C})$, $I(\mathcal{C})$ is the homogeneous ideal of $\mathcal{C}$, $\mathfrak{M} = (x_0, \dots, x_{p+1})R$ and $H^1_\mathfrak{M}(R)$ is the first local cohomology module with respect to $\mathfrak{M}$. Recall that $R$ is Cohen-Macaulay if $H^1_\mathfrak{M}(R)=0$ and is Buchsbaum if $\mathfrak{M}H^1_\mathfrak{M}(R)=0.$ More generally, for $k\geq 1,$ we say that $R$  is \textit{$k$-Buchsbaum} if $\mathfrak{M}^k H^1_\mathfrak{M}(R) = 0$ and 
\textit{strictly $k$-Buchsbaum} if $k$ is the least integer such that $\mathfrak{M}^k H^1_\mathfrak{M}(R) = 0$.


In this paper, we focus on (projective) monomial curves which serve as a basic object for several problems. These curves are defined by their generic zero $(s^d, s^{d-r_1}t^{r_1}, \dots, s^{d-r_p}t^{r_p}, t^d)$ where $r_1 < r_2 < \dots < r_p < d$ are positive integers with $\gcd(r_1, r_2, \dots, r_p, d) = 1$. On the other hand, let \( \mathcal{K}[M] \) be the subring of \( \mathcal{K}[s,t] \) generated by the set of monomials, $M = \{s^d, s^{d-r_1}t^{r_1}, \dots, s^{d-r_p}t^{r_p}, t^d\}$. If \( R \) is the homogeneous coordinate ring of the respective monomial curve, then \( R \cong \mathcal{K}[M] \). Such a curve can be seen as a projection of one-dimensional Veronese variety (see \cite{11}).  

Study of Veronese varieties and their projections dates back to Gröbner's problem, posed in \cite{19}. Gröbner showed that the Veronese variety is arithmetically Cohen-Macaulay. He posed the problem of determining conditions under which a monomial projection of Veronese variety is arithmetically Cohen-Macaulay. For simple and double projections, whose set of parameters  are obtained by deleting one or two monomials from set of all monomials of some degree $d$ in $p+2$ indeterminates, Grobner's problem was satisfactorily answered by 
Schenzel \cite{Schenzel-78} and Trung \cite{Trung-double}. Interestingly, by the classifications given in \cite{Schenzel-78} and \cite{Trung-double}, it follows that almost all such projections are arithmetically Buchsbaum.  In fact, the first non-trivial example of a non-Cohen-Macaulay ring, due to Macaulay \cite{20},  is  $k[s^4, s^3t, st^3, t^4]$ which turns out to be a Buchsbaum ring. 
It is natural to ask the conditions under which the multiple projections of Veronese variety is Buchsbaum. This was investigated in \cite{11} and \cite{Hoa-triple} for one dimensional Veronese variety and for triple projections respectively. 
The sequences in $\mathbb{P}^3$ whose  corresponding monomial curves are strictly Buchsbaum or strictly $2$-Buchsbaum are known due to Bresinsky-Schenzel-Vogel \cite{10} and  Fiorentini-Hoa \cite{22} respectively. Using the theory of affine semigroup rings, some criteria for $k$-Buchsbaumness have been found \cite{9}, \cite{11} and \cite{22}. However,
in literature, not many 
non-Cohen-Macaulay rings \( R \) are known whose corresponding monomial curves are non-smooth, see \cite{6, 9,10,11}. We recall that a curve defined by the sequence $r_1 < r_2 < \dots < r_p < d$ is smooth if and only if \( r_1 = 1 \) and \( r_p = d - 1 \).

Our results in this paper exhibit
non-smooth monomial curves in $\mathbb{P}^m$ for  $m>3$ 
 whose coordinate rings are non-Cohen-Macaulay Buchsbaum. These curves are not studied earlier. Our motivation to construct these curves is rooted in the structure of the Buchsbaum curve corresponding to $k[s^4,s^3t,st^3,t^4]$ discovered by Macaulay in $\mathbb{P}^3.$ To state our result, let $n \geq 1$, $r \geq 5$ and consider the set of monomials  
\begin{equation}\label{M-eq}
\mathcal{M} = \{s^{a_n}, s^{a_n - b_1}t^{b_1}, s^{a_n - b_2}t^{b_2}, s^{a_n - a_0}t^{a_0}, s^{a_n - a_1}t^{a_1}, \dots, t^{a_n}\}
\end{equation}
where $b_1 = r$, $b_2 = 2r - 2$ and $a_i = 2r + 2i - 1$ for all $0 \leq i \leq n$. Let $\mathcal{R}\cong\mathcal{K}[\mathcal{M}]$. 
We prove the following results.  

\begin{theorem}[Theorem \ref{XY}]\label{t1}
Suppose $5\leq r \leq 2n+5$. Then  
\begin{enumerate}[(i)]  
\item $\mathcal{R}$ is Cohen-Macaulay if $r$ is odd and $2n \geq r + 1$.  
\item $\mathcal{R}$ is non-Cohen-Macaulay Buchsbaum if $r$ is odd and $2n < r + 1$.  
\item $\mathcal{R}$ is non-Cohen-Macaulay Buchsbaum if $r$ is even.  
\end{enumerate} 
\end{theorem}
Next, we find strictly $2$-Buchsbaum non-smooth curves in $\mathbb{P}^m$ for  $m>3.$
\begin{theorem}[Corollary \ref{XZ}]\label{t2}
Suppose $n \in \mathbb{N}$ and $r = 2n + 8$. Then $\mathcal{R}$ is strictly $2$-Buchsbaum.  
\end{theorem}
In $\mathbb{P}^4$, we find strictly $k$-Buchsbaum monomial curves  for  $k\geq 1.$
\begin{theorem}[Theorem \ref{QQ}]\label{t3}
Suppose $n=1$ and $k\geq 1.$ 
Then \(\mathcal{R}\) is strictly \(k\)-Buchsbaum if and only if  
\[
r \in \{\,2 + 3k,\; 3 + 3k,\; 4 + 3k\,\}. 
\] 
\end{theorem}

It is of great interest to determine the Castelnuovo-Mumford regularity \( \operatorname{reg}(R) \), for a monomial curve $\mathcal{C}$, where $R$ is the homogeneous coordinate ring of $\mathcal{C}$. This is important because \( \operatorname{reg}(R) \) measures the complexity of the structure of graded modules. Several important invariants of graded rings and modules can be estimated via regularity. In view of the Eisenbud-Goto conjecture \cite{12}, which was proved for projective curves by Gruson-Lazarsfeld-Peskine \cite{13}, we know that \( \operatorname{reg}(R) \leq d - p \). L’vovsky gave a new bound in terms of the gap \cite{14} of the sequence $r_1,\ldots,r_p,d$. For smooth curves, Hellus-Hoa-Stückrad \cite{15} showed that \( \operatorname{reg}(R) \) is bounded by a fraction of the largest gap. In \cite{6,15,16}, one can see explicit regularity formulas for some classes of smooth and non-smooth monomial curves. In this paper, 
we establish a relation between regularity and $k$-Buchsbaumness for certain types of monomial curves and as an immediate consequence, we recover a result of \cite{22}, which states that for a smooth monomial curve, $R$ is strictly $k$-Buchsbaum if and only if $\operatorname{reg}(R) = k + 1$. In fact,   Lam-Trung \cite{6} studied some classes of non-smooth monomial curves and proved that those curves are strictly Buchsbaum if and only if $\operatorname{reg}(R)=2$. However, for the non-smooth curves defined by \eqref{M-eq}, we prove in Corollary \ref{c}   that  if $\mathcal{R} $ is strictly $k$-Buchsbaum then $\operatorname{reg}(\mathcal{R})\geq k+2$  equality holds in $\mathbb{P}^4$. The later is proved in Theorem \ref{VVV}.


Suppose $Q= (s^d, t^d)$. It is well known that  \( r_Q(R) \leq \operatorname{reg}(R) \) for any monomial curve, where \( r_Q(R) \) denote the reduction number of \( R \) with respect to \( Q  \), i.e., the least integer $m$ such that ${R}_{m+1} = Q_{m+1}$. It is of interest to find curves for which \( r_Q(R) = \operatorname{reg}(R) \), see \cite{15}. Unfortunately, this equality does not hold in general, even for smooth monomial curves \cite[Example 3.4]{15}. Whenever it does hold, it provides a combinatorial proof of the Eisenbud-Goto conjecture for monomial curves, see  \cite{17}. In Theorem \ref{C_1}, we prove that $\operatorname{reg}(\mathcal{R}) = r_Q(\mathcal{R})$, where $Q = (s^{a_n}, t^{a_n})$ where $\mathcal{R}$ is defined by the monomials in \eqref{M-eq}.  

In the study of smooth curves, one comes across Veronese variety as the Macaulayfication of the curve. Recall that
$\widetilde{R}$ is called the {\it {Macaulayfication of $R$}} if
 \ \( R \subseteq \widetilde{R}  \), \( \widetilde{R} \) is Cohen-Macaulay and \( \widetilde{R} / R \) has finite length. Although, \( R \)  has a unique finite Macaulayfication $\widetilde{R}$ for any monomial curve, 
 finding $\widetilde{R}$ is a difficult task. In \cite{1}, Goto provided a method to construct Macaulayfication for Buchsbaum rings. One can see  a numerical description of Macaulayfication in \cite{24}. Inspired by  \cite{1}, we present a generalized method for finding the Macaulayfication of an arbitrary monomial curve in Section~\ref{BE}. In particular, using this method, we construct the Macaulayfication of $\mathcal{R}$, see Remark \ref{k}, which is instrumental for our results on regularity of $\mathcal{R}$ as stated earlier. 

This paper is organised as follows. In Section \ref{section 1}, we discuss some basic facts about monomial ideals for the convenience of reader. Section \ref{AE} is dedicated to proving Theorem (\ref{t1}). In Section \ref{BE}, we study Macaulayfication of monomial curves. In Section \ref{CE}, we prove Theorem \ref{t2} along with several other results on $\operatorname{reg}(R)$ and $r_Q(R)$.  In the last section, we prove Theorem \ref{t3} and show that regularity of these curves is $k+2.$  We end the paper with examples and a question about the ring $\mathcal{R}.$

\section{Preliminaries}\label{section 1}
In this section, we gather some basic facts about homogeneous ideals and  projective closures of monomial curves. First, we recall that for an ideal $I\subseteq \mathcal{K}[x_1, \dots, x_{p+1}],$ 
 the \textit{homogenization} of $I$ is defined as the ideal
\[
I^h:= (f^h : f \in I) \subseteq \mathcal{K}[x_0, x_1, \dots, x_{p+1}],
\]
where $f^h$ denotes the homogenization of the polynomial $f$ .
Next, for an affine variety $W \subseteq \mathcal{K}^{p+1}$, the \textit{projective closure} of $W$ is the projective variety
\[
\overline{W}:= V(I_a(W)^h) \subseteq \mathbb{P}^{p+1}_\mathcal{K},
\]
where $I_a(W)^h \subseteq \mathcal{K}[x_0, \dots, x_{p+1}]$ is the homogenization of the ideal $I_a(W) \subseteq \mathcal{K}[x_1, \dots, x_{p+1}]$ of $W$. It is well known that for an ideal $I \subseteq \mathcal{K}[x_1, \dots, x_{p+1}]$, $V(I^h) \subseteq \mathbb{P}^{p+1}_\mathcal{K}$
is the projective closure of the affine variety $V_a(I) \subseteq \mathcal{K}^{p+1}$, see \cite[Chapter 8]{2}.
Let  $$\alpha'': 0=r_0<r_1<\ldots<r_{p+1}$$ be a sequence of non-negative integers such that $\gcd(r_1, \dots, r_{p+1}) = 1.$ Suppose $d=r_{p+1}$ and consider the sequence $$\alpha': d - r_0, d - r_1 \dots, d - r_{p+1}.$$ 
 Define a $\mathcal{K}$-algebra homomorphism
\[
\varphi: \mathcal{K}[x_0, x_1, \ldots, x_{p+1}] \to \mathcal{K}[s, t] \text{ by } \varphi(x_i) = s^{d - r_i} t^{r_i}, \quad 0 \le i \le p+1
\]
and let $I(\alpha)$ denote the kernel of $\varphi$. Then \[
R = \mathcal{K}[x_0, \ldots, x_{p+1}]/I(\alpha) \cong \mathcal{K}[s^d, s^{d - r_1}t^{r_1}, \ldots, s^{d - r_p}t^{r_p}, t^d].
\] Further, let $\varphi^{\prime\prime}$ and $\varphi^\prime$ denote the maps $\mathcal{K}[x_1, \ldots, x_{p+1}] \to \mathcal{K}[t] \text{ defined as } \varphi^{\prime\prime}(x_i) =  t^{r_i},~ 1 \le i \le p+1$ and  $\mathcal{K}[x_0, \ldots, x_{p}] \to \mathcal{K}[s] \text{ defined as } \varphi^{\prime}(x_i) =  s^{d-r_i}, \quad 0 \le i \le p$ respectively. We write 
$I(\alpha^{\prime\prime})$ and $I(\alpha')$ for $\ker (\varphi^{\prime\prime})$ and $\ker (\varphi^{\prime})$ respectively. Then, we have 
\begin{equation}\label{1}
\begin{aligned}
R'' = \mathcal{K}[x_1, \ldots, x_{p+1}]/I(\alpha'') \cong \mathcal{K}[t^{r_1}, t^{r_2}, \ldots, t^{r_p}, t^d] \text{ and}
\\R' = \mathcal{K}[x_0, \ldots, x_p]/I(\alpha') \cong \mathcal{K}[s^d, s^{d - r_1}, \ldots, s^{d - r_p}].
\end{aligned}
\end{equation}
Note that $I(\alpha'')$ is the vanishing ideal of the affine monomial curve $\mathcal{C}(\alpha'')$ given parametrically by $t \mapsto (t^{r_1}, \ldots, t^d)$,
and $I(\alpha')$ is the vanishing ideal of the affine curve $\mathcal{C}(\alpha')$ given parametrically by
$s \mapsto (s^d, s^{d - r_1}, \ldots, s^{d - r_p})$. Moreover,
\begin{equation}\label{2}
\begin{aligned}
I(\alpha) &= I(\alpha'')^h \quad \text{(homogenized with respect to } x_0) \\
         &= I(\alpha')^h \quad \text{(homogenized with respect to } x_{p+1}),
\end{aligned}
\end{equation}
is the vanishing ideal of the projective monomial curve $\mathcal{C}(\alpha)$ defined parametrically by $[s : t] \mapsto [s^d : s^{d - r_1}t^{r_1} : \ldots : s^{d - r_p}t^{r_p} : t^d]$. 
We end this section with the following crucial observation. 
\begin{theorem}\label{saturation}
The ideal $I = ((s^d)^m : (t^d)^\infty)$  is a primary ideal in the ring $\mathcal{K}[M]$ for all $m\geq 1.$
\end{theorem}
\begin{proof}
We need to show that $I$ has exactly one associated prime, which is equivalent to showing that all but one associated prime of the ideal $((s^d)^m)$ of $R$ contains some power of $t^d$. Let $((s^d)^m) = q_1 \cap \cdots \cap q_l$
be a minimal irredundant primary decomposition. As \(((s^d)^m)\) is a
homogeneous ideal, the ideals \(q_1,\ldots,q_l\) are homogeneous. Assume that
\(\operatorname{ht} q_i = 1\) for all \(1 \leq i \leq k\) and that \(\operatorname{ht} q_i > 1\)
for all \(k<i\leq l\). 

Saturation by \(t^d\) removes the primary components containing some power of \(t^d\). However,
if \(q_i\) contains some power of \(t^d\), then \(\operatorname{ht} q_i > 1\) since \(s^d,t^d\) is a system of
parameters. Hence $$(s^{dm}) : (t^d)^\infty = q_1 \cap \cdots \cap q_k.$$

Now we write \(p_i=\sqrt{q_i}\). These are exactly the minimal primes over \((s^{dm})\).
Now consider the integral extension \(R \subseteq K[s,t]\). For \(i=1,\ldots,k\),
let \(\mathfrak{P}_i\) be a prime ideal of \(K[s,t]\) lying over \(p_i\). We have
\(s^{dm} \in \mathfrak{P}_i\), so \(\mathfrak{P}_i=(s)\). In other words,
\(p_i=(s)\cap R\) for every \(1 \leq i \leq k\). Hence
\(q_1 \cap \cdots \cap q_k\) is primary to \((s)\cap R\).
\end{proof}
For any ring $A$ we denote by $H^i_\mathfrak{M}(A)$ the $i-$th local cohomology module of $A$ with respect to $\mathfrak{M}$ and $\ell(*)$ denote the length of a module.
\section{Buchsbaum monomial curves }\label{AE}
In this section, we will give a large class of curves which are non-smooth, non Cohen-Macaulay Buchsbaum. For a fixed positive integer $n\geq 1$, we consider the set of monomials
\[
\mathcal{M}_r^n = \{s^{a_n}, s^{a_n - b_1}t^{b_1}, s^{a_n - b_2}t^{b_2}, s^{a_n - a_0}t^{a_0}, s^{a_n - a_1}t^{a_1}, \dots, t^{a_n}\}.
\]
where $r \geq 5$ is a positive integer and \begin{equation}\label{AS}b_1 = r,~b_2 = 2r - 2 \text{ and }
a_i = 2r + 2i - 1\end{equation}for all $0\leq i\leq n.$  
For convenience, we will denote $\mathcal{M}_r^n$ simply by $\mathcal{M}$ when $r$ and $n$ are fixed and  $\mathcal{R} = \mathcal{K}[\mathcal{M}]$.
Our objective is to investigate whether $\mathcal{R}$ is a Buchsbaum ring. Note that $\mathcal{R}$ defines a non-smooth monomial curve. As discussed in \eqref{1},we have the ideals $\mathcal{I}\subset \mathcal{K}[x_0, \ldots, x_{n+3}],~\mathcal{I}'\subset \mathcal{K}[x_0, \ldots, x_{n+2}]$ and $\mathcal{I}''\subset \mathcal{K}[x_1, \ldots, x_{n+3}]$ with the following identifications:
\[
\mathcal{R} = \mathcal{K}[x_0, \ldots, x_{n+3}]/\mathcal{I} \cong \mathcal{K}[\mathcal{M}],
\]
\[
\mathcal{R}'' = \mathcal{K}[x_1, \ldots, x_{n+3}]/ \mathcal{I}'' \cong \mathcal{K}[t^{b_1}, t^{b_2}, t^{a_0}, \ldots, t^{a_n}],
\]
\[
\mathcal{R}' = \mathcal{K}[x_0, \ldots, x_{n+2}]/\mathcal{I}' \cong \mathcal{K}[s^{a_n}, s^{a_n - b_1}, s^{a_n - b_2}, s^{a_n - a_0}, \ldots, s^{a_n - a_{n-1}}].
\]
As in \eqref{2}, we have
\begin{equation*}
\begin{aligned}
\mathcal{I} &= (\mathcal{I}'')^h \quad \text{(homogenization with respect to } x_0), \\ 
           &= (\mathcal{I}')^h \quad \text{(homogenization with respect to } x_{n+3}).
\end{aligned}
\end{equation*}
\begin{theorem}\label{thm-pre5}
 In the ring $\mathcal{R}$, the ideal $I=(s^{a_n},s^{3a_n-3b_1}t^{3b_1})$ is a primary ideal.    
\end{theorem}
\begin{proof}
In view of Theorem \ref{saturation}, it is enough to prove that $I=((s^{a_n}):(t^{a_n})^\infty)$. Let $X$ be a monomial in $\mathcal{R}$ and $X.(t^{a_n})^m\in (s^{a_n})$ for some non-negative integer $m$. Without loss of generality, we may assume that
\begin{equation*}
    X=(s^{a_n})^{m_1}(s^{a_n-b_1}t^{b_1})^{m_2}(s^{a_n-b_2}t^{b_2})^{m_3}(s^{a_n-a_0}t^{a_0})^{n_0}\cdots(s^{a_{n}-a_{n-1}}t^{a_{n-1}})^{n_{n-1}}  
\end{equation*} 
where $m_1,m_2,m_3,n_i $ are nonnegative integers for $i=0,1,\dots, n-1$ and we may assume that $m_1=0$ and $m_2<3$. In this case, we want to prove that $X\in (s^{a_n}).$ We write $X^\prime$ for the monomial $X$  after the substitution $t=1.$ Thus $X' \in \mathcal{R}'$,  i.e.,   \begin{equation}\label{prime}
X^\prime =(s^{a_n-b_1})^{m_2}(s^{a_n-b_2})^{m_3}(s^{a_n-a_0})^{n_0}(s^{a_n-a_1})^{n_1}\ldots(s^{a_n-a_{n-1}})^{n_{n-1}}
\end{equation}
Since $X.(t^{a_n})^m\in (s^{a_n})$ in $\mathcal{R}$ we get that $ X^\prime \in (s^{a_n})$ in $\mathcal{R}'$. 
\begin{claim}
    There exist non negative integers $m_1',m_2',m_3',n_i',$ $i=1,...,n-1$ and $m_1'\neq 0$ such that 
\begin{eqnarray}\label{eq}
    X'=(s^{a_n})^{m_1'}(s^{a_n-b_1})^{m_2'}(s^{a_n-b_2})^{m_3'}\cdots(s^{a_n-a_{n-1}})^{n_{n-1}'}\\ 
    \mbox{ with }\label{ineq} 
     m_2+m_3+\sum_{i=0}^{n-1}n_i \geq m_1'+m_2'+m_3'+\sum_{i=0}^{n-1}n_i'.
\end{eqnarray} 
\end{claim}
\begin{proof}[Proof of Claim]
Observe that $a_n = 2r + 2n + 1,~a_n-{b_1} = r + 2n - 1,~a_{n} - b_{2} = 2n + 1,~a_n - a_i = 2n - 2i, \; i = 0, \ldots, n - 1$
and 
$a_n > a_{n}-b_1 > a_{n}-b_2 > \cdots > a_n - a_{n-1}.$   Let $$P = m_2(a_n - b_1) + m_3(a_n - b_2) + \sum_{i=0}^{n-1} n_i(a_n - a_i).$$
Then $X' = s^P$ in $\mathcal{K}[s]$. Now $X' \in (s^{a_n})\mathcal{R}' \Rightarrow X' = s^{a_n} s^{P-a_n}$ where $s^{P-a_n} \in \mathcal{R}'$. On comparing the degree of $s$ in $\mathcal{K}[s]$ we get,
 $P - a_n \geq 0.$ 
 
 Suppose $P - a_n < a_n - b_2$. Then $P-a_n$ is a non-negative integer combination of $a_n-a_i$ for $0\leq i\leq n-1$ . Notice that $a_n - a_i$, $i = 0, \ldots, n$ are all the even numbers from $0$ to $a_n - b_2$. Thus, $P - a_n$ is an even number and $P - a_n = a_n - a_j \text{ for some } j \in \{0, 1, \ldots, n\}.$
In this case, we get that $X' = s^{a_n} s^{(a_n - a_j)} \text{ and } m_1' + m_2' + m_3' + \sum_{i=0}^{n-1} n_i' = 2 \leq m_2 + m_3 + \sum_{i=0}^{n-1} n_i.$

Now suppose $P-a_n \geq a_n - b_2$.  Let $m \geq 2$ be a positive integer such that 
$(m - 1) \leq \frac{P-a_n}{a_n - b_2} < m.$ Then $ (a_n - b_2) \leq (P - a_n) - (m - 2)(a_n - b_2) < 2(a_n - b_2).$ For convenience, let
\[
L = (P - a_n) - (m - 2)(a_n - b_2).
\]
Suppose $L$ is an odd number. Then
\[
0 \leq L - (a_n - b_2) < a_n - b_2
\]
and $L - (a_n - b_2)$ is an even number which means $L - (a_n - b_2) = a_n - a_j \text{ for some } j \in \{0, \ldots, n\}.$ This gives
\[
s^P = s^{a_n + (m - 2)(a_n - b_2) + (a_n - b_2) + (a_n - a_j).}
\]
So, we take $m_1' = 1$, $m_2' = 0$, $m_3' = (m - 1)$ and $n_j' = 1$ for some $j \in \{0,1, \ldots,n\}$.
Now we show the inequality \eqref{ineq} in three cases when $m_2 = 0,1$ and $2$ respectively. First let $m_2=0$. Then 
\begin{eqnarray*}
\left( m_3 + \sum_{i=0}^{n-1} n_i \right)(a_n - b_2) \geq P& =& a_n + (m - 1)(a_n - b_2) + (a_n - a_j)
\text{ and }\\
P &>& (1 + (m - 1))(a_n - b_2) 
\end{eqnarray*}
where the last inequality follows since $ a_n > a_n - b_2 \text{ and } a_n - a_j \geq 0.$ This implies $
m_3 + \sum_{i=0}^{n-1} n_i > m
\Rightarrow m_3 + \sum_{i=0}^{n-1} n_i \geq m + 1
\Rightarrow m_3 + \sum_{i=0}^{n-1} n_i \geq m_1' + m_2' + m_3' + \sum n_i'
.$

Suppose $m_2 = 1$. Then
\begin{eqnarray*}
(a_n - b_1) + \left( m_3 + \sum_{i=0}^{n-1} n_i \right)(a_n - b_2) &\geq& P 
= a_n + (m - 1)(a_n - b_2) + (a_n - a_j) \text{ and }
\\ P &>& (a_n - b_1) + (m - 1)(a_n - b_2) 
\end{eqnarray*}
where the last inequality follows since $a_n > a_n - b_1.$ This implies $m_3 + \sum_{i=0}^{n-1} n_i > m - 1
\Rightarrow m_3 + \sum_{i=0}^{n-1} n_i \geq m
\Rightarrow m_2 + m_3 + \sum_{i=0}^{n-1} n_i \geq m_1' + m_2' + m_3' + \sum_{i=0}^{n-1} n_i'.$

Now suppose $m_2 = 2$. Then
\begin{eqnarray*}
 2(a_n - b_1) + \left( m_3 + \sum_{i=0}^{n-1} n_i \right)(a_n - b_2) &\geq& P 
= a_n + (m - 1)(a_n - b_2) + (a_n - a_j)
\\ \text{ and }  P &>& 2(a_n - b_1) + (m - 2)(a_n - b_2)
\end{eqnarray*}
where the last inequality follows since  $a_n + (a_n - b_2) > 2(a_n - b_1)$. This implies  \\$\left( m_3 + \sum_{i=0}^{n-1} n_i \right) > m - 2
\Rightarrow m_3 + \sum_{i=0}^{n-1} n_i \geq m - 1
\Rightarrow m_2 + m_3 + \sum_{i=0}^{n-1} n_i \geq m_1' + m_2' + m_3' + \sum_{i=0}^{n-1} n_i'.$

\vspace{0.5cm}
Now assume that $L$ is an even number. We have that $(a_n - b_2) - (a_n - a_0) \leq L - (a_n - a_0) < 2(a_n - b_2)-(a_n-a_0) \Rightarrow 1 \leq L - (a_n - a_0) < (a_n - b_2) + 1.$ Since $L - (a_n - a_0)$ is an even number and is less than $(a_n-b_2)$, 
we get that $
L - (a_n - a_0) = a_n - a_j
$
 for some $j \in \{0, 1, \ldots, n-1\}$. Therefore, $P - a_n - (m - 2)(a_n - b_2) - (a_n - a_0) = a_n - a_j$ which implies 
\[
s^P = s^{a_n + (m - 2)(a_n - b_2) + (a_n - a_j) + (a_n - a_0)}.
\]
So, we take $m_1' = 1,~m_2' = 0,~ m_3' = (m - 2)$. Further, if $j=0$, we choose $n_0'=2$ and if $j\neq 0$, we choose  $n_0' = 1$ and $n_j'=1$. Now we show the inequality \eqref{ineq}.  First, suppose \( m_2 = 0 \). Then,
\begin{eqnarray*}
    \left( m_3 + \sum_{i=0}^{n-1} n_i \right)(a_n - b_2) \geq P &=& a_n + (m - 2)(a_n - b_2) + (a_n - a_0) + (a_n - a_j) \text{ and }
\\P &>& 2(a_n - b_2) + (m - 2)(a_n - b_2) 
\end{eqnarray*}
where the last inequality follows since $a_n+ (a_n - a_0 )> 2(a_n - b_2) \text{ and } a_n - a_j > 0$. This implies $\left( m_3 + \sum_{i=0}^{n-1} n_i \right) > m
\Rightarrow m_3 + \sum_{i=0}^{n-1} n_i \geq m + 1=m_1' + m_2' + m_3' + \sum_{i=0}^{n-1} n_i'.$

Now suppose  \( m_2 = 1 \). Then,
\begin{eqnarray*}
    (a_n - b_1) + \left( m_3 + \sum_{i=0}^{n-1} n_i \right)(a_n - b_2)&\geq&P= a_n + (m - 2)(a_n - b_2) + (a_n - a_0) + (a_n - a_j) \\ \text{ and }
P &>& (a_n - b_1) + (a_n - b_2) + (m - 2)(a_n - b_2) 
\end{eqnarray*}
where the last inequality follows since  $a_n + (a_n - a_0) > (a_n-b_1)+(a_n - b_2)$. This implies $m_3 + \sum_{i=0}^{n-1} n_i > m - 1
\Rightarrow m_3 + \sum_{i=0}^{n-1} n_i \geq m
\Rightarrow m_2 + m_3 + \sum_{i=0}^{n-1} n_i \geq m_1' + m_2' + m_3' + \sum_{i=0}^{n-1} n_i'.$
Finally let \( m_2 = 2 \). Then
\begin{eqnarray*}
    2(a_n - b_1)+ \left(m_3 + \sum_{i=0}^{n-1} n_i\right)(a_n - b_2) &\geq& P= a_n + (m - 2)(a_n - b_2) + (a_n - a_0) + (a_n - a_j) 
\\ \text{ and } P &>& (m - 2)(a_n - b_2) + 2(a_n - b_1) 
\end{eqnarray*}
where the last inequality follows since $a_n + (a_n - a_0) > 2(a_n - b_1)$. This implies \\$\left( m_3 + \sum_{i=0}^{n-1} n_i \right) \geq m - 1
\Rightarrow m_2 + m_3 + \sum_{i=0}^{n-1} n_i \geq m_1' + m_2' + m_3' + \sum_{i=0}^{n-1} n_i'.$ This completes the proof of the claim. 

\end{proof}
Now from \eqref{eq}, we get that $X'-(s^{a_n})^{m_1'}(s^{a_n-b_1})^{m_2'}(s^{a_n-b_2})^{m_3'}(s^{a_n-a_0})^{n_0'}...(s^{a_n-a_{n-1}})^{n_{n-1}'}\\=0\implies x_1^{m_2}x_2^{m_3}\cdots x_{n+2}^{n_{n-1}}-x_0^{m_1'}x_1^{m_2'}\cdots x_{n+2}^{n_{n-1}'} \in \mathcal{I}'$. So on homogenizing this element with respect to $x_{n+3}$, we get that $X \in (s^{a_n})$ in $\mathcal{R}$. Thus  $((s^{a_n}):(t^{a_n})^\infty) \subseteq I$. Now for the other side containment, we have that $s^{a_n} \in ((s^{a_n}):(t^{a_n})^\infty)$ and $3(a_n-b_1)-a_n > a_n-b_2.$ Suppose $3(a_n-b_1)-a_n$ is an even number, say $2q$ for some positive integer $q \geq 2$. Then $s^{3(a_n-b_1)}=s^{a_n}s^{(a_n-a_{n-1})q}$. So, after homogenization, we get $s^{3(a_n-b_1)}t^{3b_1}\in ((s^{a_n}):(t^{a_n})^\infty)$ .  Suppose $3(a_n-b_1)-a_n$ is an odd number. Equivalently $3(a_n-b_1)-a_n -(a_n-b_2)$ is an even number, say $2q$ for some non-negative integer $q$ then, by similar arguments, we get that $s^{3(a_n-b_1)}t^{3b_1}\in ((s^{a_n}):(t^{a_n})^\infty)$. Therefore, $((s^{a_n}):(t^{a_n})^\infty) = I$.
\end{proof}
\begin{corollary}\label{cor-pre1}
 In the ring $\mathcal{R}$, the ideal $J=(s^{ka_n},s^{(k-1)a_n}s^{3a_n-3b_1}t^{3b_1}),~k\in \mathbb{N}$  is a primary ideal. Moreover $J= ((s^{a_n})^k:(t^{a_n})^\infty).$
\end{corollary}
\begin{proof}
    Let $X\in ((s^{a_n})^k:(t^{a_n})^\infty) $ be a monomial. Then $X(t^{a_n})^m=s^{ka_n}.L$ for some non-negative integer $m$ and $L\in \mathcal{R}$,
    $\frac{X}{s^{(k-1)a_n}}(t^{a_n})^m\in (s^{a_n})\mathcal{R} \implies X\in (s^{a_n})^k$ or $X\in s^{(k-1)a_n}s^{3a_n-3b_1}t^{3b_1}$ from Theorem \ref{thm-pre5}. Thus $((s^{a_n})^k:(t^{a_n})^\infty) \subseteq J$ and $s^{3a_n-3b_1}t^{3b_1}\in ((s^{a_n}):(t^{a_n})^\infty) \implies s^{(k-1)a_n}s^{3a_n-3b_1}t^{3b_1} \in ((s^{a_n})^k:(t^{a_n})^\infty) $. Hence 
   $((s^{a_n})^k:(t^{a_n})^\infty)=(s^{ka_n},s^{(k-1)a_n}s^{3a_n-3b_1}t^{3b_1}).$
   \end{proof}
Let $\mathfrak{m}=(s^{a_n},s^{{a_n}-b_1}t^{b_1},s^{a_n-b_2}t^{b_2},s^{a_n-a_0}t^{a_0},\dots,t^{a_n})$ be the homogeneous maximal ideal of $\mathcal{R}$.
Radical of the ideal $\mathcal{L}=(s^{2a_n},t^{a_n})$ is $\mathfrak{m}$, so,  $\mathcal{L}$ is a primary ideal. We obtain the following theorem which will be needed later.
\begin{theorem}\label{thm-pre6}
Suppose $r\leq 5+2n$. Then
\begin{enumerate}[(i)]
    \item The ideal $(s^{2a_n})$ of $\mathcal{R}$ has exactly one minimal
    associated prime, and the corresponding primary component is
    $(s^{2a_n},
            s^{a_n}s^{3a_n-3b_1}t^{3b_1})$.
    Moreover, there are choices of $r$ and $n$ with which
    $s^{a_n}s^{3a_n-3b_1}t^{3b_1}$ is not a minimal generator of this
    primary component.

    \item  If $(s^{2a_n})$ has an embedded component, we can take it
    to be $\mathcal{L}$.
\end{enumerate}

\end{theorem}
The following example illustrates the redundancy in the primary
decomposition. 
\begin{example}
    Let $r=5$ and $n=3$. Then the ideal $(s^{30})$ of
 $\mathcal{K}[s^{15},s^{10}t^5,s^7t^8,s^6t^9,s^4t^{11},s^2t^{13},t^{15}]$
is primary. This can be verified using Macaulay2.
\end{example}
 Theorem \ref{thm-pre6} follows immediately from Corollary \ref{cor-pre1} and the next proposition. 
 
\begin{proposition}\label{lemma-pre1}
    If $r\leq 5+2n$ then $s^{3a_n-3b_1}t^{3b_1}t^{a_n} \in (s^{a_n}).$
\end{proposition}
\begin{proof}
Suppose $r=2r'$ is an even number where $r'$ is a positive integer. Then 
\begin{equation*}s^{3a_n-3b_1}t^{3b_1}t^{a_n}=s^{a_n}s^{2(a_n-b_2)}t^{2b_2}s^{a_n-a_{n-(r'-2)}}t^{a_{n-(r'-2)}}\in (s^{a_n}). \end{equation*}
Now suppose $r$ is an odd number, say $r=2r'+1,~r'$ is a positive integer. If $n\geq r'+1$ then \begin{equation}\label{What}
    s^{3a_n-3b_1}t^{3b_1}t^{a_n}=
        s^{2a_n}s^{a_n-a_{r'+1}}t^{a_{r'+1}}t^{a_n}\in (s^{a_n}).
        \end{equation} 
        Now suppose $n<r'+1$. On the other hand, we have $n\geq \frac{r-5}{2}=\frac{2r'-5}{2}\implies n\geq r'-2$, i.e., $r'+1>n\geq r'-2$.  Then 
       \begin{equation*} 
            s^{3a_n-3b_1}t^{3b_1}t^{a_n}=
            \begin{cases}
        s^{a_n}s^{3(a_n-b_2)}t^{3b_2}\in (s^{a_n}) \text{ when } r'-2 = n,  \\
        s^{a_n}s^{a_n-b_2}t^{b_2}s^{2(a_n-a_0)}t^{2a_0}\in (s^{a_n}) \text{ when } r'-1 = n, \\
        s^{a_n}s^{a_n-b_2}t^{b_2}s^{a_n-a_0}t^{a_0}s^{a_n-a_1}t^{a_1}\in (s^{a_n}) \text{ when } r'= n .
        \end{cases}
        \end{equation*}
\end{proof}
Now we obtain a class of Buchsbaum rings through the next theorem.    
    \begin{theorem}\label{thm-pre7}
         Suppose $r\leq 5+2n$. Then  $\mathcal{R}=\mathcal{K}[\mathcal{M}]\subseteq \mathcal{K}[s,t]$ is a  Buchsbaum ring of $\dim$ 2. 
    \end{theorem}
    \begin{proof}
      Clearly $\mathcal{R}$ is an integral domain and  $\dim \mathcal{R}=2$. Consider the ring $A=\mathcal{R}/(s^{2a_n})$. We have $\dim A=1$ and $(0)=(s^{2a_n},s^{a_n}s^{3a_n-3b_1}t^{3b_1})/(s^{2a_n}) \cap \mathcal{L}/(s^{2a_n})$ is a primary decomposition of the ideal $(0)$ in $A$ by Theorem \ref{thm-pre6}. Now from \cite[Proposition 3.13]{4},  $H_\mathfrak{m}^0(A)=(s^{2a_n},s^{a_n}s^{3a_n-3b_1}t^{3b_1})/(s^{2a_n}).$
 Now we observe the following relations which implies that $\mathfrak{m}$ annihilates $H_\mathfrak{m}^0(A)$.
    \begin{equation}\label{q_1}
s^{a_n}s^{3a_n-3b_1}t^{3b_1}s^{a_n-b_1}t^{b_1}=s^{3a_n}s^{a_n-a_0}t^{a_0}s^{a_n-a_1}t^{a_1};
       \end{equation}
       \begin{equation}\label{q_2}
        s^{a_n}s^{3a_n-3b_1}t^{3b_1}s^{a_n-b_2}t^{b_2}=s^{2a_n}s^{a_n-b_1}t^{b_1}s^{2(a_n-a_0)}t^{2a_0};
       \end{equation}
       \begin{equation}\label{q_3}
           s^{a_n}s^{3a_n-3b_1}t^{3b_1}s^{a_n-a_0}t^{a_0}=s^{2a_n}s^{a_n-b_1}t^{b_1}s^{a_n-b_2}t^{b_2}s^{(a_n-a_1)}t^{a_1};
       \end{equation}
       \begin{equation}\label{q_4}
       s^{a_n}s^{3a_n-3b_1}t^{3b_1}s^{a_n-a_i}t^{a_i}=s^{2a_n}s^{a_n-b_1}t^{b_1}s^{a_n-b_2}t^{b_2}s^{a_n-a_{i+1}}t^{a_{i+1}} 
    \end{equation}
    for $i=1,\dots,n-1$ and $s^{a_n}s^{3a_n-3b_1}t^{3b_1}t^{a_n} \in (s^{2a_n})$ from Proposition \ref{lemma-pre1}.
 Since $\mathfrak{m}$ annihilates $H_\mathfrak{m}^0(A)$, we have  that  $A$ is a Buchsbaum ring, see \cite[Lemma 2.1]{5} which implies that   $\mathcal{R}$ is a Buchsbaum ring by \cite[Theorem 3]{3}.  
    \end{proof}
  Now we find class of non Cohen-Macaulay Buchsbaum rings which is the main result of this section. 
    \begin{theorem}\label{XY}
        Suppose $2n+5\geq r\geq 5$. Then we have the following.
        \begin{enumerate}[(i)]
            \item $\mathcal{R}$ is Cohen-Macaulay if $r$ is odd and $2n\geq r+1$.
            \item $\mathcal{R}$   is   non Cohen-Macaulay Buchsbaum if $r$ is odd and $2n <r+1$. 
             \item  $\mathcal{R}$ is non Cohen-Macaulay Buchsbaum if $r$ is an even number.
            \end{enumerate}
    \end{theorem}
    \begin{proof}
    By Theorem \ref{thm-pre7}, $\mathcal{R}$ is Buchsbaum in all cases. 
    \begin{enumerate}[(i)]
        \item Let $r$ be odd, say $r=2r'+1$ with $r'$ a positive integer and $2n\geq r+1 \implies n \geq r'+1$. Then $s^{3a_n-3b_1}t^{3b_1} \in (s^{a_n})$, see (\ref{What}), which implies that the ideal $(s^{a_n})$ is a primary ideal by Theorem \ref{thm-pre5}. Therefore, every zero divisor of $\mathcal{R}/(s^{a_n})$ is in the radical of $(s^{a_n})$. So, $t^{a_n}$ is not a zero divisor of $\mathcal{R}/(s^{a_n})$. Therefore,   $\mathcal{R}/(s^{a_n})$ is Cohen Macaulay which implies that $\mathcal{R}$ is Cohen Macaulay.
       \item \label{AC} Let $r=2r'+1$. Now suppose $2n <r+1 \iff n < r'+1$. We  prove that $s^{3a_n-3b_1}t^{3b_1} \notin(s^{a_n})$. Suppose, if possible, $s^{3a_n-3b_1}t^{3b_1} \in(s^{a_n})$. Then \begin{equation}\label{AB}
           s^{3a_n-3b_1}t^{3b_1}=(s^{a_n})^{m_1}(s^{a_n-b_1}t^{b_1})^{m_2}(s^{a_n-b_2}t^{b_2})^{m_3}(s^{a_n-a_0}t^{a_0})^{n_0}\cdots (t^{a_{n}})^{n_{n}} 
       \end{equation}
       where $m_1,m_2,m_3,n_i $ are non-negative integers for $i=0,1,\dots, n$ and $m_1\neq 0$. This gives that  $(x_1)^3-(x_0)^{m_1}(x_1)^{m_2}(x_2)^{m_3}(x_3)^{n_0}(x_4)^{n_1}\cdots (x_{n+3})^{n_{n}} \in \mathcal{I}.$ Since $(x_1)^3\notin \mathcal{I}$, we get  
$m_1+m_2+m_3+\sum_{i=0}^{n}n_i=3$. 
 
Since $a_n>a_n-{b_1}>a_{n} - b_{2} >a_n - a_0>\cdots>a_n-a_{n-1}$, on comparing the degrees of $s$ in $\mathcal{K}[s,t]$ in equation (\ref{AB}),    we get that $3(a_n-b_1)\neq a_n+(a_n-b_1)+c$ for $c\in\{a_n,a_n-{b_1} ,a_{n} - b_{2} ,a_n - a_i, i = 0, \ldots, n \}$. Moreover, we have $3(a_n-b_1)>a_n+2(a_n-b_2).$ So  the equality $3(a_n-b_1)=2a_n+c$ with $c\in \{a_n,a_n-{b_1} ,a_{n} - b_{2} ,a_n - a_i, i = 0, \ldots, n \} \implies 3(a_n-b_1)-2a_n \geq 0 \implies n \geq r'+1$ which gives a contradiction. Therefore, $s^{3a_n-3b_1}t^{3b_1} \notin(s^{a_n})$ but $s^{3a_n-3b_1}t^{3b_1}.t^{a_n} \in (s^{a_n})$ by Proposition \ref{lemma-pre1}. Hence  $(s^{a_n})$ is not a primary ideal. Now $H_\mathfrak{m}^0(\mathcal{R}/(s^{a_n}))=(s^{3a_n-3b_1}t^{3b_1}) \neq 0$ by \cite[Proposition $3.13$]{4} which implies that $\mathcal{R}/(s^{a_n})$ is not a Cohen Macaulay ring. Consequently, $\mathcal{R}$ is not a Cohen Macaulay ring.
      \item Let $r$ be an even number and suppose $s^{3a_n-3b_1}t^{3b_1} \in(s^{a_n})$. Then by similar arguments as  in part (\ref{AC}), we get the equation (\ref{AB}) with $m_1+m_2+m_3+\sum_{i=0}^{n}n_i=3$ and $m_1\neq 0$. Now again using the fact that  $a_n>a_n-{b_1}>a_{n} - b_{2} >a_n - a_0>\cdots>a_n-a_{n-1}$ and  comparing the degrees of $s$ in $\mathcal{K}[s,t]$,   we get  that $(3a_n-3b_1)-a_n$ is an even number but $2a_n-3b_1 \neq 2(a_n-b_1) $ or $2(a_n-b_2)$ or $(a_n-b_1)+(a_n-b_2)$ and also $2a_n-3b_1 > 2(a_n-a_0) $ which gives a contradiction. Thus $s^{3a_n-3b_1}t^{3b_1} \notin(s^{a_n})$ but $s^{3a_n-3b_1}t^{3b_1}.t^{a_n} \in (s^{a_n})$ by Proposition \ref{lemma-pre1} which gives that $(s^{a_n})$ is not a primary ideal. Now $H_\mathfrak{m}^0(\mathcal{R}/(s^{a_n}))=(s^{3a_n-3b_1}t^{3b_1}) \neq 0$ by \cite[Proposition $3.13$]{4} which implies $\mathcal{R}/(s^{a_n})$ is not a Cohen Macaulay ring and so $\mathcal{R}$ is not a Cohen Macaulay ring.
      \end{enumerate}
    \end{proof}
    \section{Macaulayfication of monomial curves}\label{BE}
In this section, we present a method to find the Macaulayfication of a projective monomial curve. In  \cite{1}, S. Goto defined the Macaulayfication of a Buchsbaum ring under certain conditions. We adopt similar techniques here to obtain the Macaulayfication of a projective monomial curve. Let $R$ be the coordinate ring of a projective monomial curve and let $\mathfrak{M}$ denote its unique homogeneous maximal ideal. By \cite[Theorem 1.1]{6}, we have that $R$ is a $k$-Buchsbaum domain of dimension $2$ for some non-negative integer $k$ and it admits a Macaulayfication, denoted by $\widetilde{R}$. This implies that $R$ is a generalized Cohen-Macaulay ring. Consequently, there exists a positive integer $l$ such that for every system of parameters $\alpha_1, \alpha_2$ of $R$, we have $(\alpha_{i-1}) : \alpha_i \subseteq (\alpha_{i-1}) : \mathfrak{M}^l~ \text{for } i = 1, 2, \text{ see \cite{7}.}$
Let \(A\) be a Noetherian graded ring with homogeneous maximal ideal \(\mathbf{m}\), and let \(I\) be an ideal of \(A\). Assume that \(I\) has a primary decomposition
$I=\bigcap_i Q_i(P_i),$ where each $Q_i(P_i)$ is \(P_i\)-primary. We define
\[
C_A(I)=\bigcap_{P_i\neq \mathbf{m}} Q_i(P_i).
\]
For an $x\in R$, we write $C(x)$ for $C_R((x))$ for convenience.

\begin{lemma}\label{4.1}
    If $R$ is $k$-Buchsbaum and $a\in \mathfrak{M}^k$  then $R/(a)$ is also $k$-Buchsbaum.
\end{lemma}
\begin{proof}
    Since $R$ is domain, $0 \to R \xrightarrow{a}R\to R/(a)\to 0$ is an exact sequence which induces the exact sequence $0\to H^0_{\mathfrak{M}}(R/(a))\to H^1_{\mathfrak{M}}(R)\to ...$ of local cohomology modules and    $\mathfrak{M}^kH^1_{\mathfrak{M}}(R)=0$ hence $\mathfrak{M}^kH^0_{\mathfrak{M}}(R/(a))=0$. Thus $R/(a)$ is $k$-Buchsbaum.
\end{proof}
\begin{lemma}\label{4.2}
 Let $l\in \mathbb{Z}^+$ be an integer such that \begin{equation}\label{trungs-result}
   (\alpha_{i-1}):\alpha_i\subseteq(\alpha_{i-1}):\mathfrak{M}^l \text{  for } i = 1, 2
 \end{equation}
 for every system of parameters    $\alpha_1,\alpha_2$ of $R$. Let $a,b$ be a system of parameters of $R$ with $b\in \mathfrak{M}^l$. Then
\begin{enumerate}[(i)]
    \item $C(a)=(a:b^t)=a:\mathfrak{M}^l$ for all $t\geq 1.$ 
   \item $C(a)^2=aC(a).$
\end{enumerate}
\end{lemma}
\begin{proof}
    \begin{enumerate}[(i)]
        \item \label{AD} Let $a,b$ be a system of parameters of $R$ and $b\in \mathfrak{M}^l$. Then $a,b^t$ is a system of parameters for $t\geq 1$. Now using 
        \eqref{trungs-result},
    we have 
        $$(a:b^t)\subseteq (a:\mathfrak{M}^l)\subseteq (a:b)\subseteq (a:b^t) $$
    which gives $(a:b^t)=(a:\mathfrak{M}^l)=(a:b)$  for $t\geq 1.$  
   Then we get $(0:b')=(0:(b')^t)=(0:(\mathfrak{M}')^l)$ in $R/(a)$ where $b'=b/(a)$ and $\mathfrak{M}'=\mathfrak{M}/(a)$ . On the other hand, as $C(a)/(a) = C_{R/(a)}(0)$ and by \cite[Proposition 3.13]{4}, we get that  $C_{R/(a)}(0) = H^0_{\mathfrak{M}'}(R/(a))=\mathop\bigcup\limits_{r>0}(0:(\mathfrak{M}')^r)=(0:(\mathfrak{M}')^{l_1})$ for large $l_1\in \mathbb{Z}^+$. Now $(0:(\mathfrak{M}')^{l_1}) \subseteq (0:(b')^{l_1})=(0:(\mathfrak{M}')^l) \implies  H^0_{\mathfrak{M}'}(R/(a))=(0:(\mathfrak{M}')^l) $. Thus we get $C(a)=(a:b)=a:\mathfrak{M}^l$.

\item  It suffices to show that $C(a)^2 \subseteq aC(a)$. Let $f, g \in C(a)$ and $a,b$ be a system of parameters of $R$ with $b\in \mathfrak{M}^l$. Then we may express $bf = ax$ and $bg = ay$ for some $x, y \in R$. On the other hand, we have $fg = az$ for some $z \in R$. Hence $a(b^2z) = b^2(fg) = a^2(xy)\implies b^2z = a(xy)$ as $a$ is $R$-regular. Thus $z \in (a : b^2)$ and consequently $z \in C(a)$ as $C(a) = a : b^2$ from part (\ref{AD}). Therefore $fg \in aC(a)$ and so we have $C(a)^2 \subset aC(a)$ as required.
\end{enumerate}
\end{proof}
Let \(a\) be a nonzero element of \(R\) and $Q(R)$ be the total quotient ring of $R$. It is easy to see that $R\subseteq \tilde{R}\subseteq Q(R)$ where $$\tilde{R}=\{\frac{x}{a}\in Q(R):x\in C(a)\}=a^{-1}C(a).$$ 
\begin{lemma}
    Suppose that  $\Tilde{R}$ is as above then
    \begin{enumerate}[(i)]
        \item $\Tilde{R}$ is a subring of $Q(R).$
        \item $\Tilde{R} \cong C(a)$ as $R$ modules.
    \end{enumerate}
\end{lemma}
\begin{proof}
    \begin{enumerate}[(i)]
        \item This follows from Lemma \ref{4.2}.
        \item Let $f:C(a)\to \Tilde{R}$ be the map defined by $f(x)=\frac{x}{a}.$ Then $f$ is the required isomorphism.
    \end{enumerate}
\end{proof}
\begin{theorem}\label{H}
Suppose that $\tilde{R}$ is as above. Then $$H^1_\mathfrak{M}(R) \cong \Tilde{R}/R.$$ Moreover,  $\Tilde{R}$ is the Macaulayfication of $R.$
\end{theorem}
\begin{proof}
Since \(\widetilde{R}\) and \(\widetilde{R}/R\) are \(R\)-modules, we may consider the following short exact sequence of \(R\)-modules:
\begin{equation}\label{d'}
0 \longrightarrow R \longrightarrow \widetilde{R} \longrightarrow \widetilde{R}/R \longrightarrow 0.
\end{equation}
Applying the local cohomology functor $H^i_\mathfrak{M}(-)$, we observe that $H^0_\mathfrak{M}(\Tilde{R}/R) = \Tilde{R}/R$,
since $H^0_\mathfrak{M}(\tilde{R}/R) = \mathop\bigcup\limits_{r > 0} (0 : \mathfrak{M}^r)$ and from Lemma \ref{4.2}(\ref{AD}), we have $(0 : \mathfrak{M}^l) = \Tilde{R}/R$ for some $l \in \mathbb{Z}^+$. 

We have that $R$ is a $k$-Buchsbaum domain for some $k\geq 0$. Let \(a\) be a nonzero element of $\mathfrak{M}^k$. Now consider the short exact sequence of $R$-modules:
\begin{equation}\label{a'}
0 \to (a) \to C(a) \to C(a)/(a) \to 0
\end{equation}
Since $R \cong (a)$ as $R$-modules, it follows that
\begin{equation}\label{b'}
H^i_\mathfrak{M}(R) \cong H^i_\mathfrak{M}((a)) = H^i_\mathfrak{M}(a) \quad \text{for all } i > 0.
\end{equation}

It is known that $C(a)/(a) = H^0_\mathfrak{M}(R/(a))$, we have $\dim H^0_\mathfrak{M}(R/(a)) = 0$ and $H^i_\mathfrak{M}(C(a)/(a))=0$ for all $i>0$. Since $a$ is a non-zero divisor of both $R$ and $C(a)$, we have  $H^0_\mathfrak{M}(C(a)) = 0$. Applying $H^i_\mathfrak{M}(-)$ to the exact sequence \eqref{a'} yields
\begin{equation}\label{eee}
    0 \to C(a)/(a) \to H^1_\mathfrak{M}(a) \to H^1_\mathfrak{M}(C(a)) \to 0.
\end{equation}
This implies $H^i_\mathfrak{M}(a) \cong H^i_\mathfrak{M}(C(a))\text{ for } i \geq 2$ as $R$-module.
Now consider the short exact sequence of $R$-modules:
\begin{equation}\label{a}
0 \to (a) \xrightarrow{f} R \to R/(a) \to 0,
\end{equation}
where $f:(a)\to R$ is the inclusion map, and $g:R\longrightarrow (a)$
is the multiplication map by $a$, that is, $g(x)=ax,$ for all $x\in R.$
Then $(fg)(x)=ax,$ for all $x\in R.$


From \eqref{a}, applying $H^i_\mathfrak{M}(-)$ gives
\begin{equation}\label{b}
0 \to H^0_\mathfrak{M}(R/(a)) \to H^1_\mathfrak{M}(a) \xrightarrow{H^1_{\mathfrak M}(f)} H^1_\mathfrak{M}(R) \to H^1_\mathfrak{M}(R/(a)) \to \cdots.
\end{equation}
Since \( g \) is an isomorphism,  we have the following commutative diagram:
 \[
 \begin{tikzcd}
 & H^1_\mathfrak{M}(R) \arrow[dl, "H^1_{\mathfrak M}(g)"'] \arrow[dr, "\lambda_a"] & \\
 H^1_\mathfrak{M}(a) \arrow[rr, "H^1_{\mathfrak M}(f)"] & & H^1_\mathfrak{M}(R)
 \end{tikzcd}
 \]
 We know that $a \cdot H^1_\mathfrak{M}(R) = 0$ as $a \in \mathfrak{M}^k$. Thus the map $H^1_{\mathfrak M}(f)$ is the zero map. Therefore, from \eqref{b}, we get the exact sequence 
\[
0 \to H^0_\mathfrak{M}(R/(a)) \to H^1_\mathfrak{M}(a) \to 0,
\]
and hence
\begin{equation}\label{c'}
H^0_\mathfrak{M}(R/(a)) \cong H^1_\mathfrak{M}(a) \cong H^1_\mathfrak{M}(R)
\end{equation}
as $R$-modules. Now $H^0_\mathfrak{M}(R/(a))\cong C(a)/(a)$ and from \eqref{c'}, we get that   $\ell(H^0_\mathfrak{M}(R/(a))) = \ell(H^1_\mathfrak{M}(a))$ and we know that an injective homomorphism between  two  artinian modules of the same length is an isomorphism.  
Therefore from \eqref{eee} we get $H^1_\mathfrak{M}(C(a)) = 0.$ Hence, we obtain
\[
H^0_\mathfrak{M}(C(a)) = 0, \quad H^1_\mathfrak{M}(C(a)) = 0, \quad H^i_\mathfrak{M}(C(a)) = H^i_\mathfrak{M}(R) \text{ for } i \geq 2.
\]
Since $\dim_R C(a) = \operatorname{depth}_R C(a) = 2$, it follows that $C(a)$ is a Cohen-Macaulay $R$-module. Thus, $\tilde{R}$ is a Cohen-Macaulay $R$-module. Now, applying $H^i_\mathfrak{M}(-)$ to the sequence \eqref{d'}, we get the exact sequence 
\[
0 \to \tilde{R}/R \to H^1_\mathfrak{M}(R) \to 0
\]
which gives the desired isomorphism $H^1_\mathfrak{M}(R) \cong \Tilde{R}/R.$
Now from \cite[Theorem 1.1]{6} we get that $\Tilde{R}$ is the Macaulayfication of $R$.
\end{proof}

Using our results, we obtain the Macaulayfication of  $\mathcal{R}=\mathcal{K}[\mathcal{M}].$ 
\begin{remark}\label{k}
Let $\mathcal{R}=\mathcal{K}[\mathcal{M}]$ be as in Section \ref{AE} and suppose that $s^{3a_n - 3b_1} t^{3b_1} \notin (s^{a_n})$. Then, for a positive integer $k'$, 
\[
(s^{k'a_n}) = (s^{k'a_n}, s^{(k'-1)a_n} s^{3a_n - 3b_1} t^{3b_1}) \cap (s^{k'a_n}, t^{ma_n})
\]
is a primary decomposition of the ideal $(s^{k'a_n})$
where $m \in \mathbb{Z}^+$. Consequently, $
\widetilde{\mathcal{R}} = \mathcal{K}[\mathcal{\widetilde{M}}]$ is the Macaulayfication of $\mathcal{R}$ where $\widetilde{\mathcal{M}}=\mathcal{M}\cup \{ s^{2a_n - 3b_1} t^{3b_1}\}.$
\end{remark}
\begin{proof}
From Theorem~\ref{saturation} and Corollary~\ref{cor-pre1}, we have that $(s^{k'a_n}, s^{(k'-1)a_n} s^{3a_n - 3b_1} t^{3b_1}) = (s^{k'a_n} : (t^{a_n})^\infty)$
is a primary ideal. Moreover, for  some sufficiently large $m \in \mathbb{Z}^+$, it is easy to verify that $(s^{3a_n - 3b_1} t^{3b_1}) \cdot (t^{a_n})^m \in (s^{a_n})$, see Theorem~\ref{thm-pre5}. The ideal $(s^{k'a_n}, t^{ma_n})$ is primary since its radical ideal is maximal.  Therefore, we conclude that $(s^{k'a_n}) = (s^{k'a_n}, s^{(k'-1)a_n} s^{3a_n - 3b_1} t^{3b_1}) \cap (s^{k'a_n}, t^{ma_n})$ is a primary decomposition. 

Suppose $\mathcal{R}$ is $k$-Buchsbaum for some $k\geq 1$. Then we have 
$$C(s^{ka_n} \mathcal{R}) = (s^{ka_n}, s^{(k-1)a_n} s^{3a_n - 3b_1} t^{3b_1}).$$
From Theorem \ref{H}, it follows that the Macaulayfication of $\mathcal{R}$ is
$
\widetilde{\mathcal{R}} = \mathcal{K}[\widetilde{\mathcal{M}}]$ where $ \widetilde{\mathcal{M}}=\mathcal{M}\cup \{ s^{2a_n - 3b_1} t^{3b_1}\}.$
\end{proof}
\section{Regularity of monomial curves}\label{CE}
\noindent
In this section, we discuss the reduction number and the Castelnuovo-Mumford regularity of a projective monomial curve. Let $R \cong \mathcal{K}[M]$ be the coordinate ring of a projective monomial curve, where $M = \{s^d, s^{d-r_1}t^{r_1}, \dots, s^{d-r_p}t^{r_p}, t^d\}$.  The reduction number of $R$ with respect to the ideal $Q = (s^d, t^d)$, denoted by $r_Q(R)$,  is defined as the smallest integer $m$ such that $R_{m+1} = Q_{m+1}$.   Let $\widetilde{R}$ be the Macaulayfication of $R$, given by $\widetilde{R} = \mathcal{K}[\widetilde{M}]$, where $\widetilde{M}=M \cup N$ suppose $N$ consists of monomials of the same degree, say $ld$, for some positive integer $l$. Define $R' = \mathcal{K}[N]$. We have that $R$ has standard grading with the monomials in $M$  having degree one and $\deg s^\alpha t^\beta=(\alpha+\beta)/d$ 
for all $s^\alpha t^\beta\in R$ and $R'$ has standard grading with the monomials in $N$  having degree one and $\deg s^\alpha t^\beta=(\alpha+\beta)/ld$ for all $s^\alpha t^\beta\in R'.$ 
We write $R = \bigoplus\limits_{i \geq 0} R_i$ and  $R' = \bigoplus\limits_{i \geq 0} R_i'$. 
\begin{lemma}\label{lemma reg}
Suppose $R_2' \subseteq R_{2l}$ or $R_2' \subseteq R_l R_1'$. Then $R$ is strictly $k$-Buchsbaum if and only if $a(\widetilde{R}/R) := \max\{n \in \mathbb{Z} \mid (\widetilde{R}/R)_n \neq 0\} = k + l - 1.$
\end{lemma}
\begin{proof}
We have that $\widetilde{R}$ is a graded $R$-module with grading $\widetilde{R} = \bigoplus\limits_{i \geq 0} \widetilde{R}_i$, where
\[
\widetilde{R}_i = 
\begin{cases}
R_i & \text{for all } i < l, \\
R_l \oplus R_1' & \text{for } i = l, \\
R_{l+1} \oplus R_1 R_1' & \text{for } i = l+1, \\
\vdots & \\
R_{2l} \oplus R_l R_1' \oplus R_2' & \text{for } i = 2l,\\
\vdots
\end{cases}
\]
As $R_2' \subseteq R_{2l}$ or $R_2' \subseteq R_l R_1'$, we get that for all $i \geq 0,~
\widetilde{R}_{l+i} = R_{l+i} \oplus R_i R_1'$.
 $R$ is strictly $k$-Buchsbaum if and only if $\mathfrak{M}^k H^1_{\mathfrak{M}}(R) = 0 \quad \text{but} \quad \mathfrak{M}^{k-1} H^1_{\mathfrak{M}}(R) \neq 0.$ Using Theorem \ref{H}, it is equivalent to say that $\mathfrak{M}^k \widetilde{R} \subseteq R \quad \text{but} \quad \mathfrak{M}^{k-1} \widetilde{R} \nsubseteq R$
if and only if $\widetilde{R}_{k+l} \subseteq R_{k+l} \quad \text{but} \quad \widetilde{R}_{k+l-1} \nsubseteq R_{k+l-1}.$
Equivalently, the largest integer $m$ such that $(\widetilde{R}/R)_m\neq 0$ is $k + l - 1$, i.e., $a(\widetilde{R}/R) = k + l - 1$.
\end{proof}


 Let $r_Q(\widetilde{R})$ denote the least integer $m$ such that $\widetilde{R}_{m+1} = (Q \cdot \widetilde{R})_{m+1}$ where $Q=(s^d,t^d).$
\begin{corollary}\label{corollary reg}
Suppose $R$ is strictly $k$-Buchsbaum and $R_2' \subseteq R_{2l}$ or $R_2' \subseteq R_l R_1'$. 
Then $\operatorname{reg}(R) = \max\{k + l,\, r_Q(\widetilde{R})\}.$
\end{corollary}

\begin{proof}
The result follows directly from \cite[Proposition 1.8]{6} and Lemma~\ref{lemma reg}.
\end{proof}

\begin{corollary}\label{cor-reg-1}
Suppose $R$ is a smooth monomial curve. Then $R$ is strictly $k$-Buchsbaum if and only if $\operatorname{reg}(R) = k + 1$.
\end{corollary}

\begin{proof}
For a smooth monomial curve $R$,  $\widetilde{R}$ is the Veronese subring of $\mathcal{K}[s,t]$ generated by all monomials of degree $d$, so we have $l = 1$. It is also straightforward to verify that $R_2' \subseteq R_1 R_1'$ or $R_2' \subseteq R_2$. From \cite[Example 1.9]{6}, we know that $r_Q(\widetilde{R}) = 1$. Hence, by Corollary~\ref{corollary reg}, we obtain the desired result.
\end{proof} 
\begin{corollary}\label{c}
Let $\mathcal{R}=\mathcal{K}[\mathcal{M}]$ as in Section \ref{AE} be strictly $k$-Buchsbaum . Then $\operatorname{reg}(\mathcal{R}) = \max\{k + 2,\, r_Q(\widetilde{\mathcal{R}})\}.$
\end{corollary}

\begin{proof}
By Remark \ref{k}, we have $N=\{s^{(2a_n - 3b_1)} t^{(3b_1)}\}$, i.e., $l = 2$ and $s^{2(2a_n - 3b_1)} t^{2(3b_1)} \in \mathcal{R}$, which implies $\mathcal{R}_2' \subseteq \mathcal{R}_{2l}$. Applying Lemma~\ref{lemma reg} and Corollary~\ref{corollary reg}, we obtain the result.
\end{proof}

In the next theorem, we obtain a numerical criteria for $k$-Buchsbaumness. Our result is a generalization of \cite[Theorem 1.2]{6} in which case, $l=k=1$. 
\begin{theorem}\label{theorem:5.1}
The ring $R$ is strictly $k$-Buchsbaum if and only if $kG_R + G_{\widetilde{R}} = (k + l)G_R$,
where $$G_R = \{a \in \mathbb{Z}_{\geq 0} \mid s^{d-a}t^a \in M\}, \quad G_{\widetilde{R}} = lG_R \cup \{a \in \mathbb{Z}_{\geq 0} \mid s^{ld - a}t^a \in N\}$$
and $lG_R$ denotes the set of all sums of $l$ elements from $G_R$.
\end{theorem}
\begin{proof}
Clearly $(k + l)G_R\subseteq kG_R + G_{\widetilde{R}}$. We know that $R_k \cdot \widetilde{R}_l \subseteq \widetilde{R}_{k+l}$. By Lemma~\ref{lemma reg}, $R$ is strictly $k$-Buchsbaum if and only if $\widetilde{R}_{k+l} \subseteq R_{k+l}$. Define the sets $A = \{s^{kd - a}t^a \mid a \in kG_R\}, ~B = \{s^{ld - a}t^a \mid a \in G_{\widetilde{R}}\},~C = \{s^{(k + l)d - a}t^a \mid a \in (k + l)G_R\}$.
Then $R_k$, $\widetilde{R}_l$ and $R_{k+l}$ are $\mathcal{K}$-linear spans of the monomials in $A$, $B$ and $C$, respectively. Thus, $R$ is strictly $k$-Buchsbaum if and only if the product set $A \cdot B \subseteq C$, which is equivalent to $kG_R + G_{\widetilde{R}} \subseteq (k + l)G_R$. 
\end{proof}
As a consequence of the previous theorem, we find a class of strictly $2$-Buchsbaum rings. 
\begin{corollary}\label{XZ}
Suppose $n\in\mathbb{N}$ and $r = 2n + 8$ in \eqref{AS}. Then $\mathcal{R}$ is strictly $2$-Buchsbaum.
\end{corollary}
\begin{proof}
We compute that $G_{\mathcal{R}} = \{0,\, r,\, 2r - 2,\, 2r - 1,\, 2r + 1,\, 2r + 3,\, \ldots,\, 2r + 2n - 1\}$ and from Remark~\ref{k}, we have $l = 2$ and $G_{\widetilde{\mathcal{R}}} = 2G_{\mathcal{R}} \cup \{3r\}.$ 
We show that $G_{\mathcal{R}}+G_{\mathcal{\tilde{R}}} \nsubseteq 3G_{\mathcal{R}}$. It follows from Theorem~\ref{theorem:5.1} that $\mathcal{R}$ is not $1$-Buchsbaum. Let $q = (2r + 2n - 1) + 3r = 12n + 39\in G_{\mathcal{R}}+G_{\mathcal{\tilde{R}}}$. However, observe that $3(2r - 2) > q$ and $0 + 2(2r + 2n - 1) < q$, so the only possible representation of $q$ is of the form $q = r + a + b$ for some $a, b \in G_{\mathcal{R}} \setminus \{0, r\}$. Since $q$ is an odd number, one of $a$ or $b$ must be $2r - 2$. But $(2r - 2) + r + (2r + 2n - 1) < q$, hence $q \notin 3G_{\mathcal{R}}$. 

We again use Theorem \ref{theorem:5.1} with $k=2.$ By comparing the degrees of \( t \) in equations \eqref{q_1}, \eqref{q_2}, \eqref{q_3} and \eqref{q_4}, we obtain the following inclusion, \begin{equation}\label{e1} G_{\mathcal{R}}\setminus\{2r+2n-1\}+\{3r\}\subseteq 3G_{\mathcal{R}}\implies G_{\mathcal{R}}+G_{\mathcal{R}}\setminus\{2r+2n-1\}+\{3r\}\subseteq 4G_{\mathcal{R}}.\end{equation} 
Now let $q' = 2(2r + 2n - 1) + 3r\in 2G_{\mathcal{R}}+\{3r\}$. We can write $q' = 2(2r - 2) + (2r - 1) + (2r + 2n - 2)\in 4G_{\mathcal{R}}$. Now for an arbitrary element $g_1+g_2+3r\in 2G_{\mathcal{R}}+\{3r\}$ other than $q',$  it follows from  \eqref{e1}, that $g_1+g_2+3r\in 4G_{\mathcal{R}}.$ Hence, we get that $2G_{\mathcal{R}}+\{3r\}\subseteq 4G_{\mathcal{R}}$ which implies $2G_{\mathcal{R}}+G_{\mathcal{\tilde{R}}}\subseteq 4G_{\mathcal{R}}$. Therefore, by Theorem \ref{theorem:5.1},  $\mathcal{R}$ is strictly $2$-Buchsbaum.
\end{proof}

\begin{theorem}\label{C_1}
   Let $\mathcal{R}=\mathcal{K}[\mathcal{M}]$ be as in Section \ref{AE}. Then $\operatorname{reg}(\mathcal{R}) = r_Q(\mathcal{R})$, where $Q = (s^{a_n}, t^{a_n})$.
\end{theorem}
\begin{proof}
   If $\mathcal{R}$ is Cohen–Macaulay or Buchsbaum, then the result follows from \cite{18}. We may assume that $\mathcal{R}$ is strictly $k$-Buchsbaum for $k\geq 2.$ From Remark \ref{k}, we have that  $(s^{a_n})=
   (s^{a_n}, s^{3a_n - 3b_1}t^{3b_1}) \cap (s^{a_n}, t^{ma_n})$ 
  is the primary decomposition of $(s^{a_n})$ for some $m \in \mathbb{Z}^+$. From Theorem \ref{H}, we have that 
  $\mathfrak{m}^k \widetilde{\mathcal{R}} \subseteq \mathcal{R}$ but $\mathfrak{m}^{k-1} \widetilde{\mathcal{R}} \nsubseteq \mathcal{R}$, where $\mathfrak{m}$ is as defined in Section \ref{AE}. Now, using Remark \ref{k}, we can say that 
$$\mathfrak{m}^k(s^{2a_n - 3b_1}t^{3b_1}) \subseteq \mathcal{R} \text{ and } \mathfrak{m}^{k-1}(s^{2a_n - 3b_1}t^{3b_1}) \nsubseteq \mathcal{R}.$$ This implies $s^{2a_n - 3b_1}t^{3b_1}t^{ka_n} \in \mathcal{R}$ but $s^{2a_n - 3b_1}t^{3b_1}t^{(k-1)a_n} \notin \mathcal{R}$, see equations \eqref{q_1}, \eqref{q_2}, \eqref{q_3} and \eqref{q_4}. This implies that $s^{3a_n - 3b_1}t^{3b_1}t^{ka_n} \in (s^{a_n})$ but $s^{3a_n - 3b_1}t^{3b_1}t^{(k-1)a_n} \notin (s^{a_n})$. Therefore, there exists a monomial $P \in \mathcal{R}$ such that  $$s^{a_n} \cdot P=s^{3a_n - 3b_1}t^{3b_1}t^{ka_n} \in \mathfrak{m}^{k+3}\implies P\in \mathfrak{m}^{k+2}.$$  
On the other hand, clearly $P \notin (t^{a_n})$. Now, if possible $P \in (s^{a_n})$, then on comparing the degrees of $s$, we obtain that $3(a_n - b_1) - 2a_n = 2n - (r+1)$
must be even. Otherwise, we cannot express $2n - (r+1) \leq 2n + 1$ as a sum of elements from the set $A = \{2r + 2n - 1,\, r + 2n - 1,\, 2n + 1,\, 2n,\, 2n - 2,\, 2n - 4,\, \ldots,\, 2,\, 0\}$. Now suppose there exists $z \in A$ such that $s^{3a_n - 3b_1} = s^{2a_n} \cdot s^z$ in $\mathcal{R}'$. After homogenization, this gives  $s^{3a_n - 3b_1}t^{3b_1} \in (s^{a_n})\mathcal{R}$  which implies that $\mathcal{R}$ is Cohen–Macaulay, contradicting our assumption.  So, $P \notin (s^{a_n})$.
 This gives  $P\notin  Q\mathfrak{m}^{k+1}$, however
   $P \in \mathfrak{m}^{k+2}$ which means $r_Q(\mathcal{R})\geq k+2$. It is well known that $r_Q(\Tilde{\mathcal{R}}) \leq r_Q(\mathcal{R}) \leq \operatorname{reg}(\mathcal{R})$. So, the equality $\operatorname{reg}(\mathcal{R}) = r_Q(\mathcal{R})=r_Q(\Tilde{\mathcal{R}})$ follows from Corollary~\ref{c}.
\end{proof}
\section{Application}\label{WWW}
Let \(\mathcal{R}\) be same as defined in Section \ref{AE} with $n=1$, i.e., \(\mathcal{R} = k[\mathcal{M}]\) is in $\mathbb{P}^4$. So, \[
\mathcal{M} = \{\, s^{2r+1},\ s^{r+1}t^r,\ s^3t^{2r-2},\ s^2t^{2r-1},\ t^{2r+1} \,\}.
\] In this section, we characterize the \(k\)-Buchsbaum property and the Castelnuovo-Mumford regularity of \(\mathcal{R}\)  in \(\mathbb{P}^4\) using the results from Sections \ref{BE} and \ref{CE}.

\begin{theorem}\label{QQ}
Let \(\mathcal{R} = k[\mathcal{M}]\) be in $\mathbb{P}^4$ and $k\geq1$. Then \(\mathcal{R}\) is strictly \(k\)-Buchsbaum if and only if $r \in \{\, 2 + 3k,\ 3 + 3k,\ 4 + 3k \,\}.$
\end{theorem}

\begin{proof}
By Remark \ref{k}, we obtain a Macaulayfication \(\widetilde{\mathcal{R}} = k[\widetilde{\mathcal{M}}]\) of \(\mathcal{R}\), where
$\widetilde{\mathcal{M}} = \mathcal{M} \cup \{\, s^{r+2}t^{3r} \,\}.$
From Theorem \ref{H}, we have
$H^1_{\mathfrak{m}}(\widetilde{\mathcal{R}}) = \widetilde{\mathcal{R}} / \mathcal{R}$
where \(\mathfrak{m}\) denotes the maximal homogeneous ideal of \(\mathcal{R}\).  
Thus, \(\mathcal{R}\) is strictly \(k\)-Buchsbaum if and only if $\mathfrak{m}^k \widetilde{\mathcal{R}} \subseteq \mathcal{R}$
and \(k\) is the smallest positive integer with this property.  
This means that $k$ is the smallest positive integer such that \(\mathfrak{m}^k (s^{r+2}t^{3r}) \subseteq \mathcal{R}\).

From equations \eqref{q_1}, \eqref{q_2} and \eqref{q_3}, it follows that $k$ is the smallest positive integer such that
$(t^{2r+1})^k (s^{r+2}t^{3r}) \in \mathcal{R}.$ Equivalently, 
$$(t^{2r+1})^k (s^{r+2}t^{3r}) = (s^{2r+1})^{m_1} (s^{r+1}t^r)^{m_2} (s^3t^{2r-2})^{m_3} (s^2t^{2r-1})^{m_4}$$ where $ m_i \geq 0 \text{ for} ~1\leq i\leq 4$ and not all $m_i$ are zeros.
By comparing the powers of \(s\) in \(k[s,t]\) and applying the property of homogenization, we obtain the system:
\[m_1=0; \quad m_2=0; \quad 3m_3 + 2m_4 = r + 2\quad  \text{ and } \quad m_3 + m_4 = k + 2
\]
which gives  \(m_3 = r - 2k - 2\geq0~\text{and} ~m_4=3k+4-r\geq 0\). This is equivalent to \begin{equation}\label{eqq}
    3k+4\geq r \geq 2k+2
\end{equation} 
and $k$ is the smallest integer such that (\ref{eqq}) holds, i.e. $3(k-1)+4\geq r\geq 2(k-1)+2$ must not hold. 
Equivalently, $3k+4\geq r>3k+1$
which holds if and only if $r\in \{\, 2 + 3k,\ 3 + 3k,\ 4 + 3k \,\}. $
\end{proof}

\begin{theorem}\label{VVV}
If \(\mathcal{R}\) is strictly \(k\)-Buchsbaum in $\mathbb{P}^4$, then \(\operatorname{reg}(\mathcal{R}) = k + 2\).
\end{theorem}

\begin{proof}
From Corollary \ref{c}, we have \(\operatorname{reg}(\mathcal{R}) \ge k + 2\),  
and from Theorem \ref{C_1}, we know \(\operatorname{reg}(\mathcal{R}) = r_Q(\mathcal{R})\),  
where \(Q = (s^{2r+1}, t^{2r+1})\).  
It is therefore enough to verify that
\begin{equation}\label{BBB}
    (s^2 t^{2r-1})^{m_3} (s^3 t^{2r-2})^{m_2} (s^{r+1} t^r)^{m_1} \in (s^{2r+1},\ t^{2r+1})
\end{equation}
where \(m_1 + m_2 + m_3 = k + 3\) and \(m_i\) are non-negative integers for \(i = 1, 2, 3\). This implies that $\operatorname{reg}(\mathcal{R})=r_Q(\mathcal{R})\leq k+2.$ Note that 
\[
(s^{r+1} t^r)^4 = (s^2 t^{2r-1}) (s^{2r+1})^2 (t^{2r+1})
\]
and $(s^2 t^{2r-1})^3 = (s^3 t^{2r-2})^2 s^{2r+1}.$
Hence, if \(m_1 \ge 4\) or \(m_3 \ge 3\), then \eqref{BBB} holds. Now let \(m_1 < 4\) and \(m_3 < 3\). 
There are exactly \(12\) such triples \((m_1, m_2, m_3)\):
\[
\begin{aligned}
& (0, k+3, 0),\ (0, k+2, 1),\ (0, k+1, 2),\ (1, k+2, 0),\ (1, k+1, 1),\ (1, k, 2),\\
& (2, k+1, 0),\ (2, k, 1),\ (2, k-1, 2),\ (3, k, 0),\ (3, k-1, 1),\ (3, k-2, 2).
\end{aligned}
\]

For each \((m_1, m_2, m_3)\) above, we claim that there exists  \((p_1, p_2, p_3, p_4, p_5)\) such that
\[
(s^2 t^{2r-1})^{m_3} (s^3 t^{2r-2})^{m_2} (s^{r+1} t^r)^{m_1}
= (t^{2r+1})^{p_5} (s^2 t^{2r-1})^{p_4} (s^3 t^{2r-2})^{p_3} (s^{r+1} t^r)^{p_2} (s^{2r+1})^{p_1},
\]
where \(p_i \geq 0\)  for \(i = 1, \dots, 5\) and \(p_1, p_5\) are not both zero. From Theorem \ref{QQ}, \(\mathcal{R}\) is strictly \(k\)-Buchsbaum if and only if  
\(r \in \{\, 2 + 3k,\ 3 + 3k,\ 4 + 3k \,\}\). The values of \((p_1, p_2, p_3, p_4, p_5)\) for each \((m_1, m_2, m_3)\) can be calculated as per the following table.
\renewcommand{\arraystretch}{1.5} 
\begin{center}
\begin{tabular}{|c|c|c|c|}
\hline
 & $r=2+3k$ & $r=3+3k$ & $r=4+3k$ \\ \hline
$(m_1,m_2,m_3)$ & $(p_1,p_2,p_3,p_4,p_5)$ & $(p_1,p_2,p_3,p_4,p_5)$ & $(p_1,p_2,p_3,p_4,p_5)$ \\ \hline
$(0,k+3,0)$ & $(0,1,2,0,k)$ & $(0,1,1,1,k)$&$(0,1,0,2,k)$ \\ \hline
$(0,k+2,1)$ & $(0,1,1,1,k)$ & $(0,1,0,2,k)$&$(0,1,1,0,k+1)$ \\ \hline
$(0,k+1,2)$ & $(0,1,0,2,k)$ &$(0,1,1,0,k+1)$ &$(0,1,0,1,k+1)$ \\ \hline
$(1,k+2,0)$ & $(0,2,1,0,k)$ &$(0,2,0,1,k)$ &$(1,0,0,1,k+1)$ \\ \hline
$(1,k+1,1)$ & $(0,2,0,1,k)$ &$(1,0,0,1,k+1)$ &$(0,2,0,0,k+1)$ \\ \hline
$(1,k,2)$   & $(1,0,0,1,k+1)$ & $(0,2,0,0,k+1)$&$(1,0,0,0,k+2)$ \\ \hline
$(2,k+1,0)$ & $(0,3,0,0,k)$ &$(1,1,0,0,k+1)$ & $(1,0,k,2,0)$\\ \hline
$(2,k,1)$   & $(1,1,0,0,k+1)$ &$(1,0,k+1,0,1)$ &$(1,0,k+1,0,1)$ \\ \hline
$(2,k-1,2)$ & $(1,0,k,1,1)$ & $(1,0,k,1,1)$& $(1,0,k,1,1)$\\ \hline
$(3,k,0)$   & $(1,1,k-1,2,0)$ &$(1,1,k-1,2,0)$ &$(1,1,k-1,2,0)$ \\ \hline
$(3,k-1,1)$ & $(1,1,k,0,1)$ &$(1,1,k,0,1)$ &$(1,1,k,0,1)$ \\ \hline
$(3,k-2,2)$ & $(1,1,k-1,1,1)$&$(1,1,k-1,1,1)$ &$(1,1,k-1,1,1)$ \\ \hline
\end{tabular}
\end{center}
\vspace{0.1cm}
 This completes the proof.
\end{proof}

In fact, in all the examples of $\mathcal{R}$ we compute, the regularity turns out to be $k+2$ if $\mathcal{R}$ is strictly $k$-Buchsbaum. 
\begin{example} Let $\mathcal{R}= \mathcal{K}[s^{21},\, s^{11}t^{10},\, s^{3}t^{18},\, s^{2}t^{19},\, t^{21}]$. We substitute $r=10$ and $n=1$ in Corollary \ref{XZ} to see that  $\mathcal{R}$ is strictly $2$-Buchsbaum ring in $\mathbb{P}^4$ of dimension two. Further, the Macaulayfication of $\mathcal{R}$ is $\mathcal{\widetilde{R}}=\mathcal{K}[s^{21}, s^{11}t^{10}, s^3t^{18}, s^{2}t^{19}, t^{21}, s^{12}t^{30}]$ by Remark \ref{k} and $\operatorname{reg}(\mathcal{R})=r_Q(\mathcal{R})=4$ by Theorem \ref{VVV}.
\end{example}
\begin{example}\label{n} Let 
$\mathcal{R}=\mathcal{K}[s^{13},s^8t^5, s^5t^8,s^4t^9,s^2t^{11},t^{13}]$.
 We substitute $r=5$ and $n=2$ in Theorem \ref{XY} to see that  $\mathcal{R}$ is non-Cohen-Macaulay Buchsbaum ring of dimension two. Further, the Macaulayfication of $\mathcal{R}$ is $\mathcal{\widetilde{R}}=\mathcal{K}[s^{13},s^8t^5, s^5t^8,s^4t^9,s^2t^{11},t^{13},s^{11}t^{15}]$ by Remark \ref{k} and we get $\operatorname{reg}(\mathcal{R})=r_Q(\mathcal{R})=3$ computed by Macaulay2.
\end{example}
\begin{example}
    Let $\mathcal{R}=\mathcal{K}[s^{15},s^{10}t^5, s^7t^8,s^6t^9,s^4t^{11},s^2t^{13},t^{15}]$.
 We substitute $r=5$ and $n=3$ in Theorem \ref{XY} to see that  $\mathcal{R}$ is a Cohen-Macaulay ring of dimension two. Further, we get $\operatorname{reg}(\mathcal{R})=r_Q(\mathcal{R})=3$ computed by Macaulay2.
\end{example}
Our computations lead us to the following question for the class of projective monomial curves as defined in Section \ref{AE} for which we do not have any counter example.    
\begin{question}\label{DDD}
If $\mathcal{R}$ is a strictly $k$-Buchsbaum ring for $k\geq 1$, is it true that $\operatorname{reg}(\mathcal{R})=k+2$?
\end{question}
\section{Acknowledgments}
The authors would like to express their gratitude to the anonymous referee for the meticulous reading of the manuscript and for providing  valuable comments and suggestions, particularly regarding the proof of Theorem~\ref{saturation}. The authors would also like to thank Ritobrata Mukherjee for initial discussions on the projective closure of Gorenstein semigroups.

\vspace{1cm}


\end{document}